\documentclass[draft]{article}
\usepackage{amsmath,amssymb}
\newtheorem{Theorem}{Theorem}
\newtheorem{Definition}[Theorem]{Definition}
\newtheorem{Lemma}[Theorem]{Lemma}
\newtheorem{Proposition}[Theorem]{Proposition}

\def \bnq{\begin{equation}} \def \enq{\end{equation}}
\def \bk{\bigskip}          \def \mk{\medskip}          
\def \nk{\noindent}

\def \mB{\mathbb{B}}
\def \mC{\mathbb{C}}    \def \mE{\mathbb{E}}        \def \mF{\mathbb{F}}
\def \mH{\mathbb{H}}            
\def \mR{\mathbb{R}}        \def \mS{\mathbb{S}}        \def \mZ{\mathbb{Z}}            
                      
              \def \cE{\mathcal{E}}

\def \cK{\mathcal{K}}                  \def \cM{\mathcal{M}}
       \def \cO{\mathcal{O}}       \def \WO{\widetilde{\Omega}}
\def \cP{\mathcal{P}}                      
\def \cT{\mathcal{T}}                  
\def \kg{\mathfrak{g}}

\def\char {\mbox{char\,}}   \def\deg{\mbox{deg\,}}      
\def\dim{\mbox{dim\,}}          \def \End{\mbox{End\,}}     \def\gr{\mbox{gr\,}}
\def \Hom{\mbox{Hom\,}} \def \Id{\mbox{Id\,}}               \def \Im{\mbox{Im\,}}
         \def\Index {\mbox{Index\,}} 
\def \ker{\mbox{ker\,}}         \def \Log{\mbox{Log\,}}     \def\Re{\mbox{Re\,}}
   \def\str{\mbox{str\,}}      \def\supp{\mbox{supp\,}}
\def \Tb{T^\bullet}             \def\Tr{\mbox{Tr}^G}            \def\tr{\mbox{tr\,}}

\begin{document}
\title{ Asymptotic equivariant index of Toeplitz operators and relative index of CR structures}

\author
{Louis Boutet de Monvel\\
{\small UPMC-Paris 6, Institut de Math\'ematiques de Jussieu-UMR7586}\\{\small 4, place Jussieu, F-75005, Paris, France,\
(boutet@math.jussieu.fr)}\\
Eric Leichtnam\\
{\small CNRS, Institut de Math\'ematiques de Jussieu, UMR7586}\\{\small 4, place Jussieu, F-75005, Paris, France}\\
Xiang Tang\thanks{Research partially supported by NSF Grant 0703775}\\
{\small Department of Mathematics, Washington University}\\{\small  St. Louis, MO, 63130 USA,\
(xtang@math.wustl.edu)}\\
Alan Weinstein
\thanks{Research partially supported by NSF Grant 0707137}\\
{\small Department of Mathematics, University of California}\\{\small  Berkeley, CA 94720 USA,\
(alanw@math.berkeley.edu)}}
\date{}
\maketitle

\tableofcontents

\section{Introduction}
Let $\Omega,\Omega'$ be two bounded Stein domains (or manifolds) with smooth strictly pseudoconvex boundaries $X_0,X'_0$ (these are compact contact manifolds), and $f_0$ a contact isomorphism $X_0\to X'_0$. If $\mH_0,\mH'_0$ denote the spaces of CR functions (or distributions) on $X_0,X'_0$ (boundary values of holomorphic functions), $S,S'$ the Szeg\"o projectors,\footnote{\ The definition of $S$ requires choosing a smooth positive density on $X_0$; nothing of what follows depends on this choice.} the map $E_0: u\mapsto S'(u\circ f_0^{-1}):\mH_0\to\mH'_0$ is Fredholm (it is an elliptic Toeplitz FIO). The index of $E_0$ was introduced by Epstein \cite{cE073}, who called it the relative index of the two CR structures. A formula for the index was proposed in \cite{aW97}. A special case was established in [20], and a proof of  this index formula in the general case was given by C. Epstein \cite{cE073}, based on an analysis of the situation using the ``Heisenberg-pseudodifferential calculus". In this paper we propose a simpler proof based on equivariant Toeplitz-operator calculus, which gives a straightforward view.

It is awkward to keep track of the index in the setting of Toeplitz operators on $X_0$ and $X'_0$ alone, because we are dealing with several Szeg\"o projectors, and Toeplitz operator calculus controls the range $\mH$ of a generalized Szeg\"o projector at best up to a vector space of finite rank\footnote{\ There is no index formula for a vector bundle elliptic Toeplitz operator, although there is one for matrix Toeplitz operators - a straightforward generalisation of the Atiyah-Singer formula, cf. \cite{lB79.2}.}.

To make up for this, we use the ball $\widetilde{\Omega}\subset \mC\times\Omega$ defined by $t\bar{t}<\phi$ where $t$ is the coordinate on $\mC$, $\phi$ is a smooth defining function ($\phi=0, d\phi\neq 0$ on $X_0=\partial\Omega$, $\phi>0$ inside - note that this is the opposite sign from the usual one) chosen so that $\Log\frac1\phi$ is strictly plurisubharmonic, so that the boundary $X=\partial\widetilde{\Omega}$ is strictly pseudoconvex; such a defining function always exists,  e.g. we can choose $\phi$ strictly pseudoconvex. $X$ is then a compact contact manifold with (positive) action of the circle group $U(1)$.  We will identify $X_0$ with the submanifold $\{0\}\times X_0$ of $X$.

 We perform the same construction for $\Omega'$: we will see that there exists an equivariant germ near $X_0$ of equivariant contact isomorphism $f: X\to X'$ extending $f_0$ such that $t'\circ f$ is a positive multiple of $t$, and an elliptic equivariant Toeplitz FIO $E$ extending $E_0$, associated \footnote{\ $f$ is associated to $E$  in the same manner as a canonical map is associated to a FIO.}  to the contact map $f$; the holomorphic spaces $\mH,\mH'$ split in Fourier components $\mH_k,\mH'_k$ on which the index is repeated infinitely many times. This construction has the advantage of taking into account the geometry of the two fillings $\Omega,\Omega'$, which obviously must come into the picture.

The final result can be then expressed in terms of an asymptotic version of the relative index ($G$-index) of $E$, derived from theory of M.F. Atiyah \cite{mA74}: the asymptotic index, described in \S\ref{index}, ignores finite dimensional spaces and is well defined for Toeplitz operators or Toeplitz systems; it is also preserved by suitable contact embeddings.

\mk The asymptotic equivariant trace and index are described in \S\ref{S1},\ref{S2}. The relative index formula is described and proved in \S\ref{S3}.

\section{Equivariant trace and index}\label{S1}
\subsection{Equivariant Toeplitz Operators.}

Let $G$ be a compact Lie group with Haar measure $dg$ ($\int dg=1$), $\mathfrak{g}$ its Lie algebra, and $X$ a smooth compact co-oriented contact manifold with an action of $G$: this means that $X$ is equipped with a contact form $\lambda$ (two forms define the same oriented contact structure if they are positive multiples of each other); $G$ acts smoothly on $X$ and preserves the contact structure,  i.e. for any $g$ the image $g_*\lambda$ is a positive multiple of $\lambda$;
replacing $\lambda$ by the mean $\int g_*\lambda dg$, we may suppose that it is invariant. The associated symplectic cone $\Sigma$ is the set of positive multiples of $\lambda$ in $T^*X$, a principal $\mR^+$ bundle over $X$, a half-line bundle over $X$.

We also choose an invariant measure $dx$ with smooth positive density on $X$, so $L^2$ norms are well defined. The results below will not depend on this choice.

\mk It was shown in \cite{BG81} that there always exists an invariant generalized Szeg\"o projector $S$ which is a self adjoint Fourier-integral projector whose microsupport is $\Sigma$, mimicking the classical Szeg\"o projector. $S$ extends or restricts to all Sobolev spaces; for $s\in\mR$ we will denote by $\mH^{(s)}$ the range of $S$ in the Sobolev space $H^s(X)$, and by $\mH$ the union.

\mk  A Toeplitz operator of degree $m$ on $\mH$ is an operator of the form $f\mapsto T_Qf=SQf$, where $Q$ is a pseudodifferential operator of degree $m$. Here we use pseudodifferential operators in a strict sense, i.e. in any local set of coordinates the total symbol has an asymptotic expansion $q(x,\xi)\sim\sum_{k\ge0}q_{m-k}(x,\xi)$  where $q_{m-k}$ is homogeneous of degree $m-k$ with respect to $\xi$, and the degree $m$ and $k\ge0$ are integers \footnote{\ We will occasionally use as multipliers operators of degree $m=\frac12$ (or any other complex number), with $k$ still an integer in the expansion.}. A Toeplitz operator of degree $m$ is continuous $\mH^{(s)}\to \mH^{(s-m)}$ for all $s$.  Recall that Toeplitz operators give rise to a symbolic calculus,  microlocally isomorphic to the pseudodifferential calculus, that lives on $\Sigma$ (cf. \cite{BG81}).

In particular, the infinitesimal generators of $G$ (vector fields determined by elements $\xi\in \mathfrak{g}$) define Toeplitz operators $T_\xi$ of degree $1$ on $\mH$.  An element $P$ of the universal enveloping algebra $U(\mathfrak{g})$ acts as a higher order Toeplitz operator $P_X$ (equivariant if $P$ is invariant), and the elements of $G$ act as unitary Fourier integral operators - or ``Toeplitz-FIO''.

\mk $\mH$ (with its Sobolev counterparts) splits according to the irreducible representations of $G$: $\mH=\widehat{\bigoplus}\; \mH_\alpha$.

\mk Below we will use the following extended notions: an equivariant Toeplitz bundle $\mE$ is the range of an equivariant Toeplitz projector $P$ of degree $0$ on a direct sum $\mH^N$. The symbol of $\mE$ is the range of the principal symbol of $P$; it is an equivariant vector bundle on $X$. Any equivariant vector bundle on $X$ is the symbol of an equivariant Toeplitz bundle (this also follows from \cite{BG81}).

\subsection{$G$-trace} \label{Gtrace}

The $G$-trace and $G$-index (relative index in \cite{mA74}) were introduced by M.F. Atiyah in \cite{mA74} for equivariant pseudodifferential operators on $G$-manifolds. The $G$-trace of such an operator $A$ is a distribution on $G$, describing $\tr (g\circ A)$. Here we adapt this to Toeplitz operators. Because the Toeplitz spaces $\mH$ or $\mE$ are really only defined up to a finite dimensional space, their $G$-trace or index are ultimately only defined up to a smooth function, i.e. they are distribution singularities on $G$ (distributions mod $C^\infty$); they are described below, and renamed ``asymptotic $G$-trace or index".

\bk
If $\mE,\mF$ are equivariant Toeplitz bundles, there is an obvious notion of Toeplitz (matrix) operator $P:\mE\to\mF$, and of its principal symbol $\sigma_d(P)$ (if it is of degree $d$), a homogeneous vector-bundle homomorphism $E\to F$ over $\Sigma$.  $P$ is elliptic if its symbol is invertible; it is then a Fredholm operator $\mE^s\to \mF^{s-d}$ and has an index which does not depend on $s$.

\mk If $\mE$ is an equivariant Toeplitz bundle and $P:\mE\to\mE$ is a Toeplitz operator of trace class\footnote{\ $\dim X=2n-1$. The Toeplitz  algebra is microlocally isomorphic to the algebra of pseudodifferential operators in $n$ real variables, and operators of degree $<-n$ are of trace class.} ($\deg P <-n)$, the trace function\footnote{\ We still denote by $g$ the action of a group element $g$ through a given representation; for example if we are dealing with the standard representation on functions, $gf = f\circ g^{-1}$, also denoted by  $g_*f$, $g^{*-1}f$, or $g^{-1*}f$.} $\Tr_P(g)=\tr (g\circ P)$ is well defined; it is a continuous function on $G$. It is smooth if $P$ is of degree $-\infty$ ($P\sim0$). If $P$ is equivariant, its Fourier coefficient for the representation $\alpha$ is $\frac{1}{d_\alpha} \tr P|_{\mE_\alpha}$  (with $d_\alpha$ the dimension of $\alpha$, $\mE_\alpha$ the $\alpha$-isotypic component of $\mE$).
\begin{Definition}
We denote by ${\rm char\,}\kg\subset\Sigma$ the characteristic set of the $G$-action, i.e. the closed subcone where all symbols of infinitesimal operators $T_\xi,\xi\in\mathfrak{g}$ vanish (this contains the fixed point set $\Sigma^G$). \end{Definition}

The fixed point set $X^G$ is the base of $\Sigma^G$ because $G$ is compact (there is an invariant section). The base $Z$ of $\char\kg$ is the set of points of $X$ where all Lie generators $L_\xi,\xi\in\mathfrak{g}$, are orthogonal to $\lambda$. It contains the fixed point set $X^G$. Note that $\Sigma^G$ is always a smooth symplectic cone and its base $X^G$ is a smooth contact manifold; $\char\kg$ and $Z$ may be singular.

\mk The following result is an immediate adaptation of the similar result for pseudodifferential operators in \cite{mA74}.
\begin{Proposition}
Let $P:\mE\to\mE$ be a Toeplitz operator, with $P\sim0$ near $\char\kg$ (i.e. its total symbol vanishes near $\char\kg$).
Then $\Tr_P = \tr (g\circ P)$ is well defined as a distribution on $G$. If $P$ is equivariant, we have, in distribution sense:
\bnq \label{tralpha}
\Tr_P=\sum \frac1{d_\alpha}\  ({\tr P|_{\mE_\alpha}}) \  \chi_\alpha
\enq
where $\alpha$ runs over the set of irreducible representations, $d_\alpha$ is the dimension and $\chi_\alpha$ the character.
\end{Proposition}

We have seen above that this is true if $P$ is of trace class. For the general case, let $D_G$ be a bi-invariant elliptic operator of order $m>0$ on $G$ (e.g. the Casimir of a faithful representation, with $m=2$). Since $D_G$ is in the center of $U(\kg)$, the Toeplitz operator $D_X: \mE\to\mE$ it defines is invariant, with characteristic set $\char\kg$.

If $P\sim0$ near ${\rm char\,}\kg$, we can divide it repeatedly by $D_X$ (modulo smoothing operators) and get for any $N$:
$$
P=D^N_X Q + R\quad\rm{with}\ R\sim 0 .
$$
The degree of $Q$ is $\deg P-N\deg(D_G)$, so it is of trace class if $N$ is large enough. We set $\Tr_P= D_G^N \Tr_Q+\Tr_R$: this is well defined as a distribution; the fact that this does not depend on the choice of $D_G,N,Q,R$ is immediate.

Formula \eqref{tralpha} for equivariant operators is obvious for trace class operators, and the general case follows by application of $D_X^N$ and $D_G^N$. Note that the series in the formula converges in the distribution sense, i.e. the coefficients have at most polynomial growth.

\bk Slightly more generally, let $(\mE,d)$
$$
\dots\to \mE_j\xrightarrow{d} \mE_{i+1}\to\dots
$$
be an equivariant Toeplitz complex of finite length, i.e. $\mE$ is a finite sequence $\mE_k$ of equivariant Toeplitz bundles, $d=(d_k:\mE_k\to\mE_{k+1})$ a sequence of Toeplitz operators such that $d^2=0$.
If the (degree zero) endomorphism $P = \{P_k\}$ of the complex $\mE$ is $\sim0$ near $\char\kg$, its supertrace $\Tr_P=\sum(-1)^k\Tr_{P_k}$ is well defined; it vanishes if $P=[P_1,P_2]$ is a supercommutator with one factor $\sim0$ on $\char\kg$.

\subsection{$G$ index} \label{Gindex}

Let $\mE_0,\mE_1$ be two equivariant Toeplitz bundles. An equivariant Toeplitz operator $P:\mE_0\to\mE_1$ is $G$-elliptic (relatively elliptic in \cite{mA74}) if it is elliptic on $\char\kg$, i.e. the principal symbol $\sigma(P)$, which is a homogeneous equivariant bundle homomorphism $E_0\to E_1$, is invertible on $\char\kg$.

Then there exists an equivariant $Q:\mE_1\to\mE_0$ such that $QP\sim1_{\mE_0}, PQ\sim 1_{\mE_1}$ near $\char\kg$. The $G$-index $I^G_P$ is then defined as the distribution $\Tr_{1-QP}-\Tr_{1-PQ}$.

\mk More generally,\footnote{\ This reduces to the case of a single operator where the complex is concentrated in degrees 0 and 1.}
an equivariant complex $\mE$ as above is $G$-elliptic if the principal symbol $\sigma(d)$ is exact on $\char\kg$. Then there exists an equivariant Toeplitz operator $s=(s_k:\mE_k\to\mE_{k-1})$ such that $1-[d,s]\sim 0$ near $\char\kg$ ($[d,s]=ds+sd$). The index (Euler characteristic) is the super trace $I^G_{(\mE,d)}=\str(1-[d,s])=\sum (-1)^j\Tr_{(1-[d,s])_j}$.

\bk For any irreducible representation $\alpha$, the restriction $P_\alpha: \mE_{0,\alpha}\to \mE_{1,\alpha}$ is a Fredholm operator with index $I_\alpha$, (resp. the cohomology $H^*_\alpha$ of $d\,|_{\mE_\alpha}$ is finite dimensional), and we have
$$
I^G_P=\sum \frac1{d_\alpha}I_\alpha \chi_\alpha\qquad
({\rm resp.}\  I^G_{(\mE,D)}=\sum_{j,\alpha} \frac{(-1)}{d_\alpha}^j \dim H^j_\alpha\ \chi_\alpha ).
$$
The $G$-index $I^G_A$ is obviously  invariant under  compact perturbation and deformation, so it only depends on the homotopy class of $\sigma(P)$ once $\mE_j$ has been chosen; it does depend on a choice of $\mE_j$ (on the projector that defines it, or on the Szeg\"o projector), because $\mE_j$ is determined by its symbol bundle only up to a finite dimensional space; this inconvenience is removed with the asymptotic index below.

\bk It is sometimes convenient to note an index as an infinite representation (mod finite representations) $\sum n_\alpha \chi_\alpha$. For the circle group $U(1)$, all simple representations are powers of the identity representation, denoted $J$, and all representations occurring as indices have a generating series
\begin{equation}\label{circle}
\sum_{k\in\mZ} n_k J^k \quad \mbox{(mod finite sums)}
\end{equation}
In fact the positive and negative parts of the series have a weak periodicity property: they are of the form $P_\pm(J^{\pm1})(1-J^{\pm k})^{-k}$ for a suitable polynomial $P_\pm$ and some integer $k$  (in other words they represent rational functions whose poles are roots of $1$, and whose Taylor series have integral coefficients).  \footnote{\ Something similar occurs for any compact group, cf. \cite{mA74}.}

\mk Here in our relative index problem, only very simple representations of the form $\sum_0^\infty mJ^k= m(1-J)^{-1}$ (for some integer $m$) will occur.

\section{K-theory and embedding}\label{S2}

A crucial point in the proof of the Atiyah-Singer index theorem consists in showing how one can embed an elliptic system $A$ in a simpler manifold where the index theorem is easy to prove, preserving the index and keeping track of the K-theoretic element $[A]$. The new embedded system $F_+A$ is analogous to a derived direct image (as in algebraic geometry), and the K-theoretic element $[F_+A]$ is the image of $[A]$ by the Bott homomorphism constructed out of Bott's periodicity theorem (cf.\cite{mA68}). Here we will do the same for Toeplitz operators. The direct image $F_+A$ is even somewhat more natural, as is its relation to the Bott homomorphism (\S\ref{embedding0}). The direct image for elliptic systems does not preserve the exact index, since this is not defined (because the Toeplitz space $\mH$ is at best only defined mod a space of finite rank); but it does preserve the asymptotic  equivariant index.

\subsection{A short digression on Toeplitz algebras}

We use the following notation: for distributions, $f\sim g$ means that $f-g$ is $C^\infty$; for operators, $A\sim B$ (or $A=B$ mod $C^\infty$) means that $A-B$ is of degree $-\infty$, i.e. has a smooth Schwartz-kernel. If $M$ is a manifold, $\Tb M$ denotes the cotangent bundle deprived of its zero section; it is a symplectic cone with base $S^*M=\Tb M/\mR_+$, the cotangent sphere bundle.

\mk As mentioned above, a compact contact $G$-manifold always possesses an invariant generalized Szeg\"o projector. More generally, if $M$ is a $G$ manifold, $\Sigma\subset\Tb M$ an invariant symplectic cone, there exists an associated equivariant Szeg\"o projector (cf \cite{BG81}). If $\Sigma\subset\Tb M,\Sigma'\subset\Tb M'$ and $f:\Sigma\to\Sigma'$ is an isomorphism of symplectic cones, there always exists an ``adapted FIO" $F$ which defines a Fredholm map $u \mapsto\tilde{F}u= S'(Fu):\mH\to\mH'$ and an isomorphism of the corresponding Toeplitz algebras  ($A\mapsto \tilde{F}A\tilde{F}^{-1}$, mod $C^\infty$).

One can choose $F$ equivariant if $f$ is. Indeed any adapted FIO can be defined using a global phase function $\phi$ on $\Tb(M\times M'{^{op}})$ such that\footnote{\ {\it op} in $M'{^{op}}$ refers to the change of sign in the symplectic form on $T^*M'$.}

1) $\phi$ vanishes on the graph of $f$, and $d\phi$ coincides with the Liouville form $\xi\cdot dx-\eta\cdot dy$ there;

2) $\Im \phi\gg0$, i.e. $\Im\phi>0$ outside of the graph of $f$, and the transversal hessian is $\gg0$; replacing $\phi$ by its mean gives an invariant phase; we may set $Ff(x)=\int e^{i\phi} a f(y) dy\,d\eta d\xi$  where the density $a(x,\xi,y,\eta)dy\,d\eta d\xi$ is a symbol, invariant and positive elliptic ($F$ is of Sobolev degree $\deg(a\,dy\,d\eta\, d\xi) -\frac34(n_x+n_y)$ (cf. H\"ormander \cite{lH71}), so $a$ is possibly of non integral degree if we want $F$ of degree $0$). The transfer map from $\mH$ to $\mH'$ is $S' F S$.

\bk If $M$ is a manifold and $X=S^*M$, the cotangent sphere, $X$ carries a canonical Toeplitz algebra, viz. the sheaf $\cE_{S^*M}$ of pseudo-differential operators acting on the sheaf $\mu$ of microfunctions. In general, if $X$ is a contact manifold, we will denote by $\cE_X$ (or just $\cE$) the algebra of Toeplitz operators on $X$.  It is a sheaf of algebras on $X$ acting on $\mu\mH=\mH$ mod $C^\infty$, which is a sheaf of vector spaces on $X$; the pair $(\cE_X,\mu\mH)$ is locally isomorphic to the pair of sheaves of pseudo-differential operators acting on microfunctions.  If $X$ is a $G$-contact manifold, we can choose the Szeg\"o projector invariant, so $G$ acts on $\cE_X$ and $\mu_X$.

For a general contact manifold, $\cE_X$ is well defined up to isomorphism, independently of any embedding - but no better than that. The corresponding  Szeg\"o projector (not mod $C^\infty$) is defined only  up to a compact operator (a little better than that - see below).

\subsection{Asymptotic trace and index}

The symbol bundles $E_j$ of the Toeplitz bundles $\mE_j$ only determines these up to a space of finite dimension (because, as mentioned above, both the projector defining them, and the Szeg\"o projector, are not uniquely determined by their symbols. However,  if $\mE,\mE'$ are two equivariant Toeplitz bundles with the same symbol, there exists an equivariant elliptic Toeplitz operator $U:\mE\to\mE'$ with quasi-inverse $V$ (i.e.  $VU\sim 1_\mE,UV\sim1_\mE'$). This may be used to transport equivariant Toeplitz operators from $\mE$ to $\mE'$: $P\mapsto Q=UPV$. Then if $P\sim0$ on $Z$,  $Q=UPV$ and $VUP$ have the same $G$-trace, and since $P\sim VUP$, we have $\Tr_P-\Tr_Q\in C^\infty(G)$.

\begin{Definition}
We define the asymptotic $G$-trace of $P$ as the singularity of $\Tr_P$ (i.e. $\Tr_P$ mod $C^\infty(G)$).
\end{Definition}
The asymptotic trace vanishes iff the sequence of Fourier coefficients of $\Tr_P$ is of rapid decrease, i.e. $O(c_\alpha)^{-m}$ for all $m$ where $c_\alpha$ is the eigenvalue of $D_G$ in the representation $\alpha$. This is the case if $P$ is of degree $-\infty$.
\begin{Definition}\label{asind}
We will say that a system $P$ of Toeplitz operators is $G$-elliptic (relatively elliptic in \cite{mA74}) if it is elliptic on $\char\kg$. When this is the case, the asymptotic $G$-index (or $\widetilde{I}^G_P$) is defined as  the singularity of $I^G_P$. (We will still denote it by $I^G_P$ if there is no risk of confusion.)
\end{Definition}

The asymptotic index  depends only on the homotopy class of the principal symbol $\sigma(P)$, and since it is obviously additive we get:
\begin{Theorem} The asymptotic index defines an additive map from $K^G(X-Z)$ to $C^{-\infty}(G)/C^\infty(G)$, where $Z$ is, as   above, the base of $\char\kg$.
\end{Theorem}

By the excision theorem $K^G(X-Z)$ is the same as $K^G_{X-Z}(X)$, the equivariant K-theory of $X$ with compact support in $X-Z$, i.e. the group of stable classes of triples $d(E,F,u)$ where $E,F$ are equivariant $G$-bundles on $X$, $u$ an equivariant isomorphism $E\to F$ defined near the set $Z$ (the equivalence relation is: $d(E,F,a)\sim 0$ if $a$ is stably homotopic (near $Z$) to an isomorphism on the whole of $X$). The asymptotic index is also defined for equivariant Toeplitz complexes, exact near $\char\kg$.

\mk Note that the sequence of Fourier coefficients $\frac{\tr  P_\alpha}{d_\alpha}$ is in any case of polynomial growth with respect  to the eigenvalues of $D$ or $D_X$; if $P\sim 0$, it is of rapid decrease. The coefficients $\frac{I_\alpha}{d_\alpha}$ of the asymptotic index are integers, so they are completely determined, except for a finite number of them, by the asymptotic index.

\mk\nk{\bf Remark}: if $V$ is a finite dimensional representation of $G$ and $V\otimes P$ or $V\otimes d$ is defined in the obvious way, we have $I^G_{V\otimes P} = \chi_V I^G_P$ (i.e. $\Index(V\otimes P)_\alpha = (V\otimes\Index P)_\alpha$, except at a finite number of places).

\mk E.g. Let $G=SU_2$ acting on the sphere $X$ of $V=\mC^2$ in the usual manner, and $E=S^m V$ the $m$-th symmetric power . Then $E\times X$ is a $G$ bundle with the action $g(v,x)=(gv,gx)$.  The CR structure on the sphere  gives rise to a first Szeg\"o projector $S_1(v\cdot f)=v \cdot S(f)$, where $S$ is the canonical Szeg\"o projector on holomorphic functions. On the other hand since $X$ is a free orbit of $G$, the bundle $E\times X$ is isomorphic to the trivial bundle $E_0\times X$ where $E_0$ is some fiber (i.e. the vector space of homogeneous polynomials of degree $m$, with trivial action of $G$). This gives rise to a second Szeg\"o projector $S_0$, not equal to the first, but giving the same asymptotic index; we recover the fact that $S^mV\otimes\sum S^kV\sim (m+1)\sum S^kV$ (= in degree $\ge m$).

\subsection{$\cE$-modules}

For the sequel, it will be convenient to use the language of $\cE$-modules. In the $C^\infty$ category, $\cE$ is not coherent; general $\cE$-module theory is therefore not practical and not usefully related to topological K-theory. We will just stick to the two useful cases below (elliptic complexes or ``good" modules).\footnote{\ Things work better in the analytic category.}.
Note also that the notion of ellipticity is slightly ambiguous; more precisely: a system of Toeplitz operators (or pseudo-differential operators) is obviously invertible mod $C^\infty$ if its principal symbol is, but the converse is not true. The notion of ``good" system below partly compensates for this; it is in fact indispensable for a good relation between elliptic systems and K-theory.

If $\cM$ is an $\cE$-module (resp. a complex of $\cE$ modules), it corresponds to the system of pseudo-differential (resp. Toeplitz) operators whose sheaf of solutions is $\Hom(\cM,\mu\mH)$; e.g.  a locally free complex of $(L,d)$ of $\cE$-modules defines the Toeplitz complex $(\mE,D)=\Hom(L,\mH)$.

More generally we will say that an $\cE$-module $\cM$ is ``good" if it is finitely generated, equipped with a filtration $\cM=\bigcup \cM_k$   (i.e. $\cE_p\cM_q=\cM_{p+q}$, $\bigcap \cM_k=0$) such that the symbol $\gr \cM$ has a finite locally free resolution. We write $\sigma(\cM)=\cM_0/\cM_{-1}$, which is a sheaf of $C^\infty$ modules on the basis $X$; since there exist global elliptic sections of $\cE$, $\gr \cM$ is completely determined by the symbol, as is the resolution.

\mk A resolution of $\sigma(\cM)$ lifts to a ``good resolution" of $\cM$, i.e. a finite locally free resolution \footnote{\ The   converse is not true: if $d$ is a locally free resolution of $\cM$, its symbol is not necessarily a resolution of the symbol of $\cM$ --  if only because filtrations must be defined to define the symbol and can be modified rather arbitrarily.} of $\cM$.

\mk It is standard that two resolutions of $\sigma(\cM)$ are homotopic, and if $\sigma(\cM)$ has locally finite locally free resolutions it also has a global one (because we are working in the $C^\infty$ category on a compact manifold or cone with compact support, and dispose of partitions of unity); this lifts to a global good resolution of $\cM$.

\bk If $\cM$ is ``good", it defines a K-theoretic element  $[\cM]\in K_Y(X)$ (where $Y$ is the support of $\sigma(\cM)$), viz. the K-theoretic element defined by the symbol of any good resolution (this does not depend on the resolution since any two such are homotopic).

\bk All this works just as well in presence of a $G$-action (if the filtration etc. is invariant).

As above (\S\ref{Gtrace}), the asymptotic $G$-trace $\Tr_A$ [using subscripts as before] is well defined if $A$ is an endomorphism of a good locally free complex of Toeplitz modules.  The same holds for a good module $\cM$: the asymptotic trace of $A\in\End\cE(M)$ vanishing near $\char\kg$ is the asymptotic trace of any lifting of $A$ to a good resolution of $\cM$. (Such a lifting, vanishing near $\char\kg$ exists and is unique up to homotopy, i.e. modulo supercommutators.) Likewise, the asymptotic $G$-index of a locally free complex exact on $Z$, or of a good $\cE$-module with support outside of $Z$, is defined: it is the asymptotic $G$-trace of the identity.

\mk Definition \ref{asind} of the asymptotic index (or Euler characteristic) extends in an obvious manner to good complexes of locally free $\cE$-modules or to good $\cE$-modules. The asymptotic $G$-index of such an object, when it is $G$-elliptic, depends only on the K-theoretic element which it defines on the base.

\mk Let us note that the asymptotic trace and index are still well defined for locally free complexes or modules with a locally free resolution, not necessarily good; in that case, what no longer works is the relation to topological K-theory on the base.

\subsection{Embedding}\label{embedding0}

If $M$ is a manifold, $\Sigma\subset\Tb M$ a symplectic subcone, the Toeplitz space $\mH$ is the space of solutions of a pseudodifferential system mimicking ${\bar\partial}_b$. If $I\subset\cE$ is the ideal generated by these operators (mod $C^\infty$), and $\cM=\cE/I$, we have $\mu\mH=\Hom_\cE(\cM,\mu)$ (as a sheaf: $f\in\Hom(\cM,\mu)\mapsto f(1)$; here as above $\mu$ denotes the sheaf of microfunctions). E.g. in the holomorphic situation, $I$ is the ideal generated by the components of $\bar{\partial}_b$.

We have $\End_\cE(\cM)=[I:I]$, the set of pseudo-differential operators $a$ such that $Ia\subset I$, acting on the right: if $a\in[I:I]$, the corresponding endomorphism of $M$ takes $f$ (mod $I$) to $fa$ (mod $I$); this vanishes iff $a\in I$. The map which takes $a\in [I:I]$ to the endomorphism $f\mapsto af$ of $\mH$ defines an isomorphism from $\End_{\cE}(\cM)$ to the algebra of Toeplitz operators mod $C^\infty$. $\cM$ is thus an $\cE_{\Tb M}-\cE_\Sigma$ bimodule (where $\cE_{\Sigma}\simeq \End_\cM$ denotes the sheaf of Toeplitz operators mod $C^\infty$).

\mk This extends immediately to the case where $\Tb M$ is replaced by an arbitrary symplectic cone\footnote{\ We use a double prime here because, eventually, we will be embedding {\em two} cones in a third one.} $\Sigma''$.  The small Toeplitz sheaf $\mu\mH$ can be realized as $\Hom_{\cE''}(\cM,\mu\mH'')$, where $\cM=\cE''/I$ and $I\subset\cE''$ is the annihilator of the Szeg\"o projector $S$ of $\Sigma$ (i.e. the null-sheaf of $I$ in $\Hom_{\cE''}(\cM,\mH'')=\mu\mH$). If $\cP$ is a (good) $\cE$-module, the transferred module is $\cM\otimes_{\cE}\cP$, which has the same solution sheaf ($\Hom(\cM\otimes\cP,\mH'')=\Hom(\cP,\Hom(\cM,\mH''))$ and $\Hom(\cM,\mH'')=\mH)$. Thus the transfer preserves traces and indices.

\bk The module $\cM=\cE''/I$ is generated by $1$ (mod $I$) and has a natural filtration, which is a good filtration: in the holomorphic case, the good resolution is dual to the complex $\bar{\partial}_b$ on $(0,*)$ forms.

In general it always has a good locally free resolution, well defined up to homotopy equivalence. In a small tubular neighborhood of $\Sigma$ one can choose this so that its symbol is the Koszul complex on $\bigwedge N'$, where $N'$ is the dual of the normal tangent bundle of $\Sigma$ equipped with a positive complex structure (as in the holomorphic case). The corresponding K-theoretic element $[\cM]\in K^G_\Sigma(\Sigma'')$ is precisely the element used to define the Bott isomorphism (with support $Y\subset\Sigma$) $K^G_Y(\Sigma)\to K^G_Y(\Sigma'')$. (Here, $Y$ some set containing  the support of $\sigma(\cM)$ and the map is the product map: $[E]\mapsto [\cM][E]$, where the virtual bundle $[E]$ on $\Sigma$ is extended arbitrarily to some neighborhood of $\Sigma$ in $\Sigma''$.)

For example if $\Sigma''$ is $\mC^N\setminus \{0\}$, with Liouville form $\Im\bar{z}\cdot dz$ and base the unit sphere  $X''=\mS^{2N-1}$), $\mH''$ is the space of holomorphic function boundary values, $\Sigma \subset \Sigma''$ consists of the nonzero vectors in the subspace $z_1=\dots=z_k=0$, and $X \subset X''$ is the corresponding subsphere, then $\mH$ consists of the functions independent of $z_1,\dots,z_k$, and $I$  is the ideal spanned by the Toeplitz operators $T_{\partial_1},\ldots T_{\partial_k}$. In this example the ideal $I$ is generated by $\bar{z}_1,\dots,\bar{z}_k$, or by $T_{\bar{z}_j},j=1\dots k$ (On the sphere we have $T_{\partial_j}=(A+N) T_{\bar{z}_j}$ with $A=T_{\sum_1^N z_j\partial_j}$). The $\cE$-module $\cM$ itself has a global resolution with symbol  the Koszul complex constructed on $\bar{z}_1,\dots,\bar{z}_k$.

\bk What precedes works exactly as well in the presence of a compact group action. If  $\cP$ is a good module with support outside of $Z$ (or a complex with symbol exact on $Z$), the transferred module has the same property ($Z\subset Z''$), and it has the same $G$-index (the $G$-index of the complex $\Hom_{\cE}(\cM,\mH) \simeq \Hom_{\cE''}(\cM'',\mH'')$).

\mk If $X,X''$ are (compact) contact $G$-manifolds, $f:X\to X''$ an equivariant embedding, $\cP$ a good $G-\cE$-module with support outside of $Z$ (the base of $\char\kg$ in $\Sigma$), or a Toeplitz complex, exact on $Z$, the transferred module on $X$ is $f_+\cP = \cM\otimes_{f_*\cE'}f_*\cP'$. This is exact outside of $f(\Sigma)$ and has the same $G$-index as $\cP$; its K-theoretic invariant $[\cP]$ is the image of $[\cP]$ by the equivariant Bott homomorphism. The K-theoretic element $[f_+\cP]\in K^G_{X-Z}(X)$ is the image of $[\cP]$ by the Bott homomorphism (it is well defined since $f(Z)\subset Z''$). Thus
\begin{Theorem} Let $f:X\to X''$ be an equivariant embedding. The Bott homomor\-phism $K^G_{X-Z}(X)\to K^G_{X''-Z''}(X'')$  commutes with  the asymptotic $G$ index.
\footnote{\ As mentioned above the interplay between the Bott isomorphism and embeddings of systems of differential or  pseudodifferential operators lies at the root of  Atiyah-Singer's proof of the  index theorem; it is described in Atiyah's works \cite{mA65,mA68,mA68.1,mA74}, cf also \cite{lB90.1} in the context of holomorphic $D$-modules, close to the Toeplitz context.}
\end{Theorem}

\bk  It is always possible to embed a compact contact  manifold in a canonical contact sphere with linear G-action. In fact, it is easier to work with the corresponding cones, as follows:
\begin{Proposition}\label{embed}
Let $\Sigma$ be a $G$-cone (with compact base), $\lambda$ a horizontal 1-form, homogeneous of degree 1, i.e. $\rho\lrcorner\lambda=0$ and $L_\rho\lambda=\lambda$, where $\rho$ is the radial vector field, generating homotheties.   Then there exists a homogeneous embedding $x\mapsto Z(x)$ of $\Sigma$ in a unitary representation space $V^c$ of $G$ such that $\lambda=\Im \bar{z}\cdot dz$.
\end{Proposition}
In the lemma, $Z$ must be homogeneous of degree $\frac12$. This applies of course if $\Sigma$ is a symplectic cone, $\lambda$ its Liouville form.   (The symplectic form is $\omega=d\lambda$ and $\lambda=\rho\lrcorner\omega$).

\mk
We first choose a smooth equivariant function $y=(y_j)$, homogeneous of degree $\frac12$, realizing an equivariant embedding of $\Sigma$ in $V-\{0\}$, where $V$ is a real unitary $G$-vector space (this always exists if the base is compact; (the coordinates $z_j$ on $V$ are homogeneous of degree $\frac12$ so that the canonical form $\Im \bar{z}\cdot dz$ is of degree $1$)). Then there exists a smooth function $x=(x_j)$ homogeneous of degree $\frac12$ such that $\lambda=2x\cdot dy$. We can suppose $x$ equivariant, replacing it  by its G-mean if need be. Since $y$ is of degree $\frac12$ we have $2\rho\lrcorner dy =y$ hence $x\cdot y=\rho\lrcorner\lambda =0$. Finally we get
$$
\lambda=\Im \bar{z}\cdot dz\quad {\rm with }\ z=x+iy .
$$

\section{Relative index}\label{S3}

 As indicated in the introduction, we are considering the index of the Fredholm map $E_0: u\mapsto S'(u\circ f_0^{-1})$ from $\mH_0$ to $\mH'_0$, where $X_0,X'_0$ are the boundaries of two smooth strictly pseudoconvex Stein manifolds $\Omega,\Omega'$, $\mH,\mH'$ the spaces CR distributions ($\ker\bar{\partial}_b$, equal to space of boundary values of holomorphic functions), $S,S'$ the Szeg\"o projectors, and $f_0$ a contact isomorphism $X_0\to X'_0$.

As announced we modify the problem and move to the larger boundaries $X,X'$ of ``balls " $|t|^2<\phi,|t'|^2<\phi'$ in $\mC\times \Omega,\mC\times \Omega'$, on which the circle group acts ($t\mapsto e^{i\lambda}t$) (\S\ref{holo}). We will see (\S\ref{collar}) that the Toeplitz FIO $E_0$ defines almost canonically an equivariant extension $F$ which is $U(1)$-elliptic, and $\Index(F|_{\mH_k})=\Index(E_0)$ for all $k$ ($\mH_k\subset \mH(X)$ is the subspace of functions $f=t^kg(x)$)
, so that our relative index $\Index(E_0)$ appears as an asymptotic equivariant index, easier to handle in the framework of Toeplitz operators.

 In \S\ref{embedding} we will show that the whole situation can be embedded in a large sphere, with action of $U(1)$ as in the examples above. In the final result (section \ref{index}) the relative index appears as the asymptotic index of an equivariant $U(1)$-elliptic Toeplitz complex on this  large sphere.  In general the equivariant index (asymptotic or not) is rather complicated to compute, but in our case the $U(1)$-action is quite simple \footnote{\ it is free on the support of the K-theoretic symbol of our complex.}, it reduces naturally to the standard Atiyah-Singer K-theoretic index formula on a symplectic ball. The result is better stated in terms of K-theory anyway, but it can be translated via the Chern character in terms of cohomology or integrals. We give here a (rather clumsy) cohomological-integral translation, essentially equivalent to the result conjectured in \cite{aW97}.

 We will also see below (\S\ref{collar}) that $f_0$ has an almost canonical extension $f$ near the boundary, well defined up to isotopy, not holomorphic but symplectic. We can then define a space $Y$ by gluing together $Y_+,Y_-$ by means of $f$. $Y$ is not a Hausdorf manifold, but it is symplectic and both $Y_+,Y_-$ carry orientations which agree on their intersection (as do the symplectic structures). We can further choose differential forms $\nu_\pm$ representatives of the Todd classes of $Y_\pm$ so that they are equal near the boundary $X_0$ (the symplectic structures agree, not the complex structures, but they define the same Todd classes).
\begin{Theorem}\label{thindex}
The relative index (index of $E_0$) is the integral $\int_Y (\nu_+ - \nu_-)$, where $\nu_\pm$ are representatives of $Todd(Y_{\pm})$ as above, so that the difference has compact support in $Y-X_0$.
\end{Theorem}
This will be explained in more detail below (\S\ref{index}). This formula is related to the Atiyah-Singer index formula on the glued space $Y$, but is not quite the same since $Y$ is not a symplectic manifold.

\bk To prove the index theorem we will give an equivalent equivariant description of the situation, where the index of $E_0$ is repeated infinitely many times, and embed everything in a large sphere where the index is given by the K-theoretic index character (\S \ref{index}).

\subsection{Holomorphic setting}\label{holo}

Let $\Omega$ be a strictly pseudoconvex domain (or Stein manifold), with smooth boundary $X_0$ ($\bar{\Omega}=\Omega\cup X_0$ is assumed to be compact);  we write $\widetilde{\Omega}\subset \mC\times\bar{\Omega}$ the ball $|t|^2<\phi$, where $\phi$ is a defining function ($\phi=0,d\phi\neq0$ on $X_0$, $\phi>0$ inside), chosen so that the boundary $X=\partial\widetilde{\Omega}$ is strictly pseudoconvex, i.e. $\Log\frac1\phi$ is strictly plurisubharmonic (i.e.
$\Im\bar{\partial}\partial\frac1\phi\gg0$).

The circle group $U(1)$ acts on $X$ by $(t,x)\mapsto (e^{i\lambda}t,x)$. We choose as volume element on $X$ the density  $d\theta\,dv$ where $dv$ is a smooth positive density on $\overline{\Omega}$ ($t=\,e^{i\theta}|t|$): this is a smooth positive density on $X$; it is invariant by the action of $U(1)$, so as the Szeg\"o projector $S$ and its range $\mH$, the space of boundary values of holomorphic functions.

The infinitesimal generator of the action of $U(1)$ is $\partial_\theta$, and we denote by $D$ the restriction to $\mH$ of
$\frac1i\partial_\theta$, which is a self-adjoint, $\ge0$, Toeplitz operator. $D$ is the restriction of $T_tT_{\partial_t}$.

\bk The model case is the sphere $\mS^{2N+1}\subset\mC^{N+1}$ with the action
$$
(t=z_0,z=(z_1,\dots,z_N)) \mapsto (e^{i\theta}t,z).
$$

\bk The Fourier decomposition of $\mH$
$$
\mH=\widehat{\oplus}_{k\ge0}\; \mH_k\qquad (\mH_k=\ker (D-k) \ )
$$
corresponds to the Taylor expansion of holomorphic functions: the $k$-th component of $f=\sum f_k(x)t^k\in\mH$ is $f_kt^k$.

$\mH_0$ identifies with the set of holomorphic functions on $X_0$ (it is the set of boundary values of holomorphic functions on $\Omega$ with moderate growth at the boundary, i.e. $|f|\le cst\ d(\cdot,X_0)^{-N}$ for some $N$, where  $d(\cdot,X_0)$ is the distance to the boundary).

\mk\nk {\bf Remark: } If $f = t^k g(x)$ with, in particular if $f\in\mH_k$, its $L^2(X)$ norm is
$$
\|f\|_{L^2(X)} = \frac{\pi}{k+1} \int_{\Omega} \phi^{k+1} |g(x)|^2 dv
$$
where as above $dv$ is the chosen smooth volume element on $\Omega$. The restriction of the Szeg\"o projector to functions of the form $t^k g(x)$ is thus identified with the orthogonal projector on holomorphic functions in $L^2(\Omega, \phi^{k+1}dv)$. Such sequences of projectors were considered by F.A. Berezin \cite{fB76} and further exploited by M. Englis \cite{mE02,mE07}, whose presentation is closely related to the one used here.

\bk For the sequel, it will be convenient to modify the factorisation $D=t\partial_t$. We begin with the easy following result.
\begin{Lemma}
Let $D=PQ$ be any factorisation where $P,Q$ are Toeplitz operators and $[D,P]=P$. Then there exists a (unique) invariant invertible Toeplitz operator $U$ such that $P=tU,Q=U^{-1}\partial_t$.
\end{Lemma}

Indeed  it is immediate that any homogeneous function $a$ on $\sigma$ such that $\frac1i\partial_\theta a=\pm a$ is a multiple $mt$ of $t$ (resp. of $\bar{t}$), with $m$ invariant. For the same reason (or by successive approximations) a Toeplitz operator $A$ such that $[D,A]=\pm A$ is a multiple of $T_tM$ (or $M'T_t$) $T_t$ with $M$ or $M'$ invariant (resp. of $T_{\partial_t}$, on the right or on the left) . Thus in the lemma above we have $P=T_tU, Q=U'T_{\partial_t}$,
where $U,U'$ are Toeplitz operators which necessarily commute with $D$, and are elliptic and inverse of each other.

Note that $D=PQ,[D,P]=P$ is equivalent to $D=PQ,[Q,P]=1$.

\mk  In particular, since $D=D^*=T_{\partial_t}^* T_t^*$, there exists a Toeplitz operator $A$ such that $T_{\partial_t}=AT_t^*$. $A$ is elliptic of degree 1 (in fact invertible), positive since $D=T_tAT_t^*$ is self-adjoint $\ge0$; it is also invariant: $[D,A]=0$.
\begin{Definition}\label{aa} We will set $\cT=T_tA^{\frac12}$; its symbol is denoted by $\sigma(\cT)=\tau$. \end{Definition}
Note that $\tau$ is homogeneous of degree $\frac12$, and $\cT$ is of degree $\frac12$, so it is not a Toeplitz operator in our strict sense, but for multiplications and automorphisms $P\mapsto UPU^{-1}$ it is just as good. We have
\bnq \cT^*=A^{\frac12}T_t^*, \ [D,\cT]=\cT. \ D=\cT\cT^* \enq
(for any other such factorisation $D=BB^*$ with $[D,B]=B$, $B$ is of degree $\frac12$, and we have $B=\cT U$ with $U$ invariant and unitary. $\cT$ is the unique Toeplitz operator giving such a factorisation and such that $\cT=T_tA'$ with $A'$ a Toeplitz operator of degree $\frac12$, $A'\gg0$).

\mk In what precedes, all $=$ signs can be replaced by $\sim$ ($=$ mod $C^\infty$); we then get local statements.

The symbol $\tau=\sigma(\cT)$ is the unique homogeneous function of degree $\frac 12$ such that $\sigma(D)=|\tau|^2, \partial_\theta \tau=i\tau, \frac \tau t>0$.

\bk We also have the following (easy) local result:
\begin{Lemma} \label{collar0}
Given any Toeplitz operator $\cK$ (mod $C^\infty$) on $\mH$ such that $D\sim \cK\cK^*,[D,\cK]=\cK$ near the boundary, there exists a unique unitary equivariant Toeplitz FIO $F$ such that $F|_{\mH_0}\sim \Id$, $F\cT \sim \cK F$.
\end{Lemma}
The geometric counterpart is: given any function $k$ on $\Sigma$ homogeneous of degree $\frac 12$ such that $\sigma(D)=k\bar{k}$ there exists a unique germ of homogeneous symplectic isomorphism $f$ such that $f|_{\Sigma_0}=\Id$, $k\circ f=\tau$. This is immediate because the two hamiltonian pairs $H_\tau, H_{\bar{\tau}}$, $H_k,H_{\bar{k}}$ define real 2-dimensional foliations, and an isomorphism $\Sigma\sim \Sigma_0\times \mC$ near $\Sigma_0$. Note that this would not work if we replaced $k,\bar{k}$ by two functions $a,b$ such that $\sigma(D)=ab, \partial_\theta a=ia$ but not $b=\bar{a}$, because then the 'foliation' defined by the Hamiltonian vector fields $H_a,H_b$, although it is formally integrable, is not real.

 The operator statement follows, e.g. by successive approximations. Note that  $F$ is completely determined by its restriction $F_0$ if it commutes with $\cT$. (In fact in $\cE_\Sigma$, the commutator sheaf of $\cT$ and $\cT^*$ identifies with the inverse image of $\cE_{\Sigma_0}$ - at least as far as the leaves of the Hamiltonian fields $H_\cT, H_{\cT^*}$ define a fibration over $\Sigma_0$: $\cE_\Sigma$ is the (completed) tensor product of the Toeplitz algebra $\mbox{Toepl}(\cT,\cT^*)$ generated by $\cT$ and $\cT^*$ and this commutator: $\cE_\Sigma\sim\cE_{\Sigma_0}\otimes\mbox{Toepl}(\cT,\cT^*)$ (in a neighborhood of $\Sigma_0$). In this statement, $(\cT,\cT^*)$ cannot be replaced by $(T_t,T_{\partial_t})$ whose commutator sheaf is only defined in the algebra of jets of infinite order along $\Sigma_0$, because the Hamiltonian leaves are complex, no longer real.)

Note that, in our case, the base of $\char\kg$ is the boundary $X_0$ (the set of fixed points), outside of which $D$ is elliptic.

\subsection{Collar isomorphisms}\label{collar}

Let now $\Omega'$ be another strictly pseudoconvex domain (or Stein manifold) with smooth boundary $X'$. We do the similar constructions $\tilde{\Omega}'$, $\mH'$, and $D'$, $\cdots$ as in the previous subsection. Let $f_0: X_0\to X'_0$ be a contact isomorphism.

We define the Fourier Toeplitz $E_0: u\mapsto S'(u\circ f_0^{-1}): \mH\to\mH'$, which is a Fredholm operator. It will be convenient to replace $E_0$ by $F_0=(E_0E_0^*)^{-\frac12}E_0$, which has the same index and is $\sim$ unitary ($E_0E_0^*$ is an elliptic $\ge0$ Toeplitz operator on $\mH'$; $(E_0E_0^*)^{-\frac12}$ is defined to be $0$ on $\ker E_0^*$  (mod $C^\infty$ would be quite enough). As for $\WO$, we construct a Toeplitz operator $\cT'$ such that $D'=\cT'\cT{'^*}, [D'\cT']=\cT', T_t^{-1}\cT'\gg0$.

Exactly as in Lemma \ref {collar}, there exists a unique (unitary) Toeplitz FIO $F$, defined near the boundary $X_0$ and mod $C^\infty$, elliptic, such that $F|_{\mH_0}=F_0$, and $F\cT\sim \cT'F$ near the boundary (mod $C^\infty$).

\mk The geometric counterpart is: there exists a unique equivariant germ of contact isomorphism $f: X\to X'$ (defined and invertible near the boundary) such that $f|_{X_0}=f_0,\ \tau'=\tau\circ f$.

We may extend $F$, using an invariant cut off Toeplitz operator, so that it vanishes (mod $C^\infty$) away from the boundary. There is an invariant FIO parametrix $F'$, i.e. $F'F\sim 1_{\mH}, FF'\sim 1_{\mH'}$, near the boundary.
\begin{Proposition} For any $k$, $F_k=F|_{\mH_k}$ has an index, equal to $\rm{Index}\, F_0$. \end{Proposition}
Proof: both $F'F$ and $FF'$ are invertible on the boundary, so have a $G$-index; the index of $F_k=F|_{\mH_k}$ is $\tr(1-F'F)_k-\tr(1-FF')_k$. Now $\cT$, resp. $\cT'$ is an isomorphism $\mH_k\to\mH_{k+1}$, resp.  $\mH'_k\to\mH'_{k+1}$, and we have $\Index(F_{k+1}A)=\Index(A'F_k)$, so  $\Index F_{k+1}=\Index F_k$, i.e. the index does not depend on $k$ and is equal to $\Index E_0$.\footnote{\ For a more general situation where $P$ is a Toeplitz operator elliptic on $X_0$, or where the canonical Szeg\"o projector is replaced by some other general equivariant one, we would only get that the index $\Index(P_k)$ is constant for $k\gg0$. Here the fact that $\Index P_k=\Index P_0$ is obvious but important.}

\mk The asymptotic index is stable by embedding; here the index is constant, and the asymptotic index of $E$ (which is essentially a Toeplitz invariant) gives the index of $F_0$ itself.

\subsection{Embedding}\label{embedding}
\begin{Theorem} Let $f:X\to X'$ be a collar isomorphism defined in some invariant neighborhood of $X_0$ in $X$. Then for large $N$ there exists equivariant contact embeddings $U:X\to \mS^{2N+1}, U':X'\to \mS^{2N+1}$ such that $U=U'\circ f$ near the boundary, and $t_X,t'_{X'}$ map to positive multiples of $t_{\mS^{2N+1}}$ (as above the contact sphere $\mS^{2N+1}$ is equipped with the $U(1)$-action $(t,z)\mapsto (e^{i\theta} t,z)$).\end{Theorem}
As usual, it will be more comfortable to work with the symplectic cones. The symplectic cone of $X$ is $\Sigma=\mR_+\times X$, where we choose the radial coordinate invariant.

 The symbol of $D$ is $\bar{\tau}\,\tau$ with $\tau/t>0$ as in Definition \ref{aa}. The Liouville form is $\Im(\bar{\tau}d\tau)+\lambda_0$ where $\lambda_0$ is a horizontal form, i.e. the pull-back of a form on $\Sigma_b=U(1)\backslash\Sigma\simeq\mR_+ \times\bar{\Omega}$ (equivalently:  $\partial_{\theta}\lrcorner \lambda_0 = L_{\partial_\theta}\lambda_0=0$).

\mk Lemma \ref{embed} provides an embedding $x\mapsto z_b(x)$ of $\Sigma_b$ in $\mC^{N'}-\{0\}$ (with the trivial action of $U(1)$). We now choose $\psi_1,\psi_2$ invariant, homogeneous of degree $0$, such that $\psi_1^2+\psi_2^2=1$, with $\supp \psi_1$ is contained in the domain of definition of $f$ and $\psi_2$ vanishing near the boundary, and we construct a new embedding $z$ in 3 pieces: $z=(z_1,z_2,z_3)$ with $z_1=\psi_1 z_0, z_2= \psi_2 z_0, z_3=0$ in $\mC^{N''}$, $N''$ to be defined below.

Since $\Im\bar{z}_jz_j\psi_jd\psi_j=0$ ($\bar{z}_jz_j\psi_jd\psi_j$ is real) we still have $\Im(\bar{z}_1\cdot dz_1+\bar{z_2}\cdot dz_2)=(\psi_1^2+\psi_2^2)\Im\bar{z_0}\cdot dz_0=\Im\bar{z_0}\cdot dz_0$ inducing $\lambda_0$. The first embedding $U=(\tau,v):\Sigma\to\mC^{1+N}$ ($N=2N'+N''$).

Similarly there exists an embedding $x'\mapsto z'_0(x')$ of  $\Sigma'_b$ in $\mC^{N''}-\{0\}$ (with the trivial action of $U(1)$). 

\mk We replace this by $z'=(z'_1,z'_2,z'_3)$ with $z'_1=\psi'_1 z_1\circ f^{-1}, z'_2=0, z'_3= \psi'_3 z'_0$ where $\psi'_1,\psi'_3$ again are invariant, homogeneous of degree $0$, $\psi{'}_1^2+\psi{'}_3^2=1$, and $\supp \psi'_1$ is contained in the domain of definition of $f^{-1}$, $\psi'_3$ vanishes near the boundary. This also defines an embedding $U'=(a',z'):\Sigma'\to\mC^{N+1}$; we have $U=U'\circ f$ near the boundary since $\psi_2,\psi'_3$ vanish there.

\subsection{Index}\label{index}

\bk We are now reduced to the case where both $U(1)$-manifolds $X,X'$ sit in a large sphere  $\mS=\mS^{2N+1}$ and coincide near the set of fixed points $\mS_0$.

As in the preceding section \ref{embedding} we can embed the $U(1)$ sheaves $\mu\mH_X,\mu\mH_{X'}$ as sheaves of solutions of two good equivariant $\cE_\mS$- modules $\cM_X,\cM_{X'}$, and the identification $F$ gives an equivariant Toeplitz isomorphism $\widetilde{F}$ near $X_0$ (we can make the construction so that $\cM_X=\cM_{X'},\widetilde{F}=\Id$ near $X_0$).

The asymptotic index then only depends on the difference element
$$
d([\cM_X],[\cM_{X'}],\sigma(\widetilde{F}))\in K^{U(1)}(\mS-\mS_0).
$$

\mk  Now $U(1)$ acts freely on $\mS-\mS_0$, with quotient space $U(1)\backslash (\mS-\mS_0)$ the open unit ball $\mB_{2N}\subset \mC^N$.  We have
\begin{Proposition}
The pull back map is an isomorphism $K(\mB)\to K^{U(1)}(\mS-\mS_0)$.

We have $K(\mB)\sim\mZ$, with generator the symbol of the Koszul complex $k_x$ at the origin (or any point of the interior), whose index is $1$.

Its pull-back is the generator of $K^{U(1)}_{\mS-\mS_0}(\mS)$: the symbol is the same, but now acting on $\mH(\mS)$. Its index is $\sum_0^\infty J^k = (1-J)^{-1}$, where (as in \eqref{circle}) $J$ is the tautological character of $U(1)$: $J(e^{i\lambda})=e^{i\lambda}$.
\end{Proposition}
The first assertion is immediate  (cf. \cite{mA74}): if $G$ is a compact group acting freely on a space $Y$, the pull back defines an equivalence from the category of vector bundles on $G\backslash Y$ to that of $G$-vector bundles on $Y$ (an inverse equivalence is given by $E\mapsto G\backslash E$), and this induces a bijection on K-theory (with supports).

The fact that $k_x$ defines the generator of $K(\mB)(=K_0(\mB)$ is just a restatement Bott's periodicity theorem. Its pullback is then the generator of $K^{U(1)}(\mS-\mS_0)$: the corresponding complex of Toeplitz operators is then the standard Koszul complex, acting on holomorphic functions, whose index is the space of holomorphic functions of $z_0=t$ alone.

Thus if $[u]\in K^{U(1)}(\mS-\mS_0)$, its asymptotic index is $m(1-J)^{-1}$ (with the notation of \S\ref{Gindex}), where the integer $m$ is the value of the K-theoretic character $K(\mB)$ on the element $[u_\mB]$ whose pull-back is $[u]$.

\mk Let us now come back to our index problem: we have constructed the difference bundle $d([\cM_X],[\cM_{X'}],\sigma(\widetilde{F}))$. We may replace $\cM_X,\cM_{X'}$ by good resolutions in small equivariant tubular neighborhoods of $X$, resp. $X'$, whose K-theoretic symbol is the Bott element -  the Koszul complex for a positive complex structure on the normal symplectic bundle of $X$, resp. $X'$. $\widetilde{F}$ lifts to the resolutions (uniquely up to homotopy), and the symbol of the lifting $u$ is an isomorphism near $X_0$ (we can make the construction so that $u=\Id$ near $X_0$), so our K-theoretic element is $[u]=d(\beta_X,\beta_{X'},u)$ (the equivariant K-theoretic element attached to the double complex defined by $u$).
\begin{Theorem} \label{thindexbis} Let $m$ be the index of $E_0$ we are investigating. Then, notations and embeddings being as above,

1) the asymptotic index of our equivariant extension $\widetilde{F}$) is the asymptotic index of the difference element $[u]=d(\beta_X,\beta_X', u)\in K^{U(1)}(\mS-\mS_0)$, where $u$ is the symbol of $\widetilde{F}$ (i.e. the identity map near $\mS_0$, where $X$ and $X'$ coincide).

2) the index $m$ itself is the value of the index character of $K(\mB)$ on the element $[u_\mB]=d(\beta_\Omega,\beta_{\Omega'}, \bar{u})$.
\end{Theorem}
The first part has just been proved. The asymptotic index is $\sim m(1-J)^{-1}$ for some integer $m$.

To prove the second we go down to $\mB_{2N}$. The bases of $X,X'$ are the embeddings $Y_+,Y_-$ of $\Omega,\Omega'$ in $\mB$, which coincide near the boundary, and as above the pullback is an isomorphism $K_{Y_\pm}(\mB)\to K^G_{X_\pm}(\mS-\mS_0)$. The Bott complexes $\beta_{X_\pm}$ descend as Bott elements $\beta_{Y_\pm}$ on $\mB$, realized as Koszul complexes of positive complex structure of the normal symplectic bundle \footnote{\ note that $Y_\pm$ are symplectic submanifolds, not complex; but all positive complex structures are homotopic.}; $u$ descends as an isomorphism near the boundary.

The index $m$ we are looking for is the K-theoretic index character of the difference element $d(\beta_{Y_+},\beta_{Y_-},u)$. This can be as usual translated in terms of cohomology, or as an integral:
$$
m = \int_\mB \omega
$$
where $\omega$ is a differential form with compact support, representative of the Chern character of our difference element $d(\beta_{Y_+},\beta_{Y_-},u)$.

We can push this down further. The construction can be made so that $u=\Id$ near the boundary, choose differential forms $\omega_\pm$ with support in small tubular neighborhoods of $Y_\pm$ so that they coincide near the boundary (so as the tubular neighborhoods), so that $\omega$ is the difference $\omega_+-\omega_-$.

The integral $\nu_\pm$ of $\omega_\pm$ over the fibers of the respective tubular neighborhoods is then a representative of the Todd class of $Y_\pm$; $\nu_+$ and $\nu_-$ coincide near the boundary, so that the difference $\nu_+-\nu_-$ has compact support in $Y=Y_+\cup Y_-$.

Finally our index $m$ is the integral $\int_Y\nu$ as announced in Theorem \ref{thindex}.

The integral can also be thought of as the constant limit $\int_{Y_{+,\epsilon}} \nu_+ -  \int_{Y_{-,\epsilon}}\nu_-$, where the subscript $\epsilon$ means that we have deleted the neighborhood $\phi<\epsilon$ in $Y_+$ and the corresponding image in $Y_-$.

\subsection{Appendix}

In this section we show how various symplectic extensions of $f_0$ are related. It is a little intriguing that, although in our proof, the extension $f$ must be chosen rather carefully so that the asymptotic index of the corresponding Toeplitz FIO $E$ is (asymptotically) the index of $E_0$, the final result, expressed as an integral on the bases glued together by means of $f$ near their boundaries, depends only on the isotopy class of $f$, which is unique.

\subsubsection{Contact isomorphisms and base symplectomorphisms}

Let $X$ be as above, with $X_0$ the fixed point set of codimension 2. Near the boundary, $X$ is identified with $X=X_0\times\mC$ and the base $U(1)\backslash X\sim\Omega$ identifies with $X_0\times \mR_+$; we have $\phi=t\bar{t}$ and the $\mC$-coordinate is $t=\sqrt{\phi}\,e^{i\theta}$ (it is smooth on $X$). The contact form is $\lambda_X= \Im(\bar{t}dt-\partial\phi)=\phi\,  d\theta+\lambda_\Omega$, where $\lambda_\Omega=-\Im\partial\phi$ is a smooth basic form. The connection form is $\gamma=d\theta -\frac{\lambda_\Omega}\phi$, and the base $\Omega=X_0\times\mR_+$ is equipped with the (basic) symplectic curvature form
$$
\mu = d\gamma \qquad(\mbox{with }\ \gamma = \frac{\lambda_\Omega}\phi, \quad \lambda_\Omega=-\Im \partial \phi)\, .
$$

We will still use the symplectic cone of $X$: this is $\Sigma=\char\kg \simeq \mR_+\times X$, with Liouville form $a\lambda_X$ and symplectic form its derivative, with the $\mR_+$ coordinate $a$ defined below: with the notation of Lemma \ref{aa}, we have $a=\sigma(A)$, i.e. $\sigma(D)=a\phi= \tau\bar{\tau}, \tau=t\sqrt{a}$ (as above $D=\frac1i T_{\partial_\theta}$ denotes the infinitesimal generator of rotations). We will also write in polar coordinates $\tau=\rho\,e^{i\theta}$ ($\rho = \sqrt{\phi\,a}$).

\mk Let $F$ be a homogeneous equivariant symplectic transformation of $\Sigma$: then $F$ preserves $\sigma(D)=\tau\bar{\tau}$, so we have necessarily $F_*\tau=u\;\tau$,  with $u$ invariant, $|u|=1$. $F$ is then completely determined by its  restriction to the boundary, since it commutes with the two real commuting hamiltonian vector fields $\Re H_\tau,\Im H_\tau$, which are linearly independent and transversal  to $\Sigma_0$.

Thus there is a one to one correspondence between unitary functions on the base $\Omega$ and germs near $\Sigma_0=\char\kg$ of equivariant symplectomorphisms inducing $\Id$ on $\char\kg$ - or equivalently of contact automorphisms of $X$ inducing $\Id$ on $X_0$.

If $F$ is such a contact automorphism, the base map $F_\Omega$ is obviously a diffeomorphism of $\Omega$ which induces $\Id$ on the boundary $X_0$ and preserves the symplectic form $\mu$.

\mk The converse is not true. If $F_\Omega$ is a smooth symplectomorphism of $\Omega$ inducing  the identity on $X_0$, we have $F_\Omega^*(\frac{\lambda_\Omega}\phi) = \frac{\lambda_\Omega}\phi + \alpha$ with $\alpha$ a closed form. It is elementary that $\alpha=c \frac{d\phi}\phi+\beta$ where $c$ is a constant and $\beta$ is smooth on the boundary. Locally on $X_0$, $F_\Omega$ lifts to $X$ or $\Sigma$: the lifting is $F:(x,\tau)\mapsto (x', \tau'=\tau e^{i\psi})$ ($\theta'=\theta+\psi$) where $\psi$ is a primitive of $\alpha$ (this is not smooth at the boundary, only continuous). It is immediate that conversely any $\alpha$ of the form above gives rise to such a contact isomorphism with smooth base map.
(on $\Sigma$ the horizontal (invariant) coordinates satisfy $H_{\tau e^{i\psi}}f=0$; the horizontal part of the Hamiltonian $H_{\tau e^{i\psi}}$ is $-i\tau e^{i\psi}(\partial_\rho -H^0_\psi)$ (with $H^0_\psi=\psi_{\xi_j}\partial_{x_j}-\psi_{x_j}\partial_{\xi_j}$); finally $\partial_\rho -H^0_\psi$ is smooth so the horizontal coordinates $(x',\xi)$ are determined by smooth differential equations.) Summing up:
\begin{Theorem} The map which to a germ of contact isomorphism $F$ (near $X_0$) assigns the invariant unitary smooth function $u$ such that $F^*\tau=\tau u$ is one to one (and continuous). In particular the homotopy class of $F$ is determined by that of $u$ (an element of $H^1(X,\mZ)$).

The map which to a smooth germ of symplectomorphism $F_\Omega$ (near $X_0$) assigns the closed one-form $\alpha = c\,\frac{d\phi}\phi+\rm{smooth}$ is one to one, the group of such symplectomorphisms is contractible. The contact lifting (which exists locally, and globally if $\alpha$ is exact) is smooth on $X_0$ iff $c=0$.
\end{Theorem}
The fact that this group is contractible (connected) simplifies the final result, namely: in the proof of Theorem \ref{thindexbis} it was essential that the base map $F_\Omega$ have a smooth symplectic extension preserving $\tau>0$; for Theorem \ref{thindex} however any symplectic $F_\Omega$ will do since these are all isotopic.

\subsubsection{Example}
(A smooth symplectic automorphism of the base does not lift to a smooth equivariant contact automorphism of the sphere.)

\bk Let $\mS$ be the unit sphere in $\mC^{N+1}$, with coordinates $x_0=t,x_1,\dots,x_N$.

$U(1)$ acts by $t\mapsto e^{i\theta}t$. The base is $\mB=\mS/U(1)$, the unit ball of $\mC^N$.

The contact form is $\Im \bar{t}dt+\lambda=\phi d\theta+\lambda$ with $\lambda=\sum\bar{x}_jdx_j$, $\phi=\bar{t}t = 1-\bar{x}x$.

The connection form is $\gamma=d\theta +\frac \lambda\phi$, its curvature is the symplectic form $\mu=d\frac\lambda\phi$ (on the interior of $\mB$).

\bk
Let $F_B$ be the diffeomorphism of $\mB$ defined by $x\mapsto x'=F_B(x)=e^{ci\phi}x$, $c$ a constant; this preserves $\phi$ and the inverse is $x=e^{-ci\phi}x'$. We have
$$
F_B^*\lambda = \Im( \bar{x} (dx+cix\,d\phi))= \lambda + c (1-\phi)d\phi
$$

Since $d(1-\phi)\frac{d\phi}\phi = 0$, $F_B$ is symplectic ($F_B^*\mu=\mu$).

But $F_B$ does not lift to a smooth equivariant contact automorphism of $\mS$: such a lifting $F$ must preserve the connection form, so it is of the form
$$
t\mapsto e^{-i\alpha}t \quad (\theta\mapsto \theta-\alpha)\quad\rm{with}\  \alpha=c\Log\phi-\phi+cst
$$
($d\alpha = c(1-\phi)\frac{d\phi}\phi$), and this is not smooth at the boundary $t=0$ if $c\neq0$.

\bk Of course the reverse works: if $F$ is a smooth equivariant contact automorphism of the sphere $\mS$ (or a germ of such near the fixed diameter $\mS_0$), the base map $F_B$ is a smooth symplectomorphism of the ball $\mB$ (up to the boundary).

\subsection{Final remarks}
1) The preceding construction applies in particular to the following situation: let $V,W$ be two compact manifolds, and $f_0$ a contact isomorphism $S^*V\to S^*W$.

We may suppose $V$ real analytic; then $S^*V$ is contact isomorphic to the boundary of small tubular neighborhoods of $V$ in its complexification. For example if $V$ is equipped  with an analytic Riemannian metric, and $(x,v)\mapsto e_x(v)$ denotes the geodesic  exponential map, the map $(x,v)\mapsto e_x(iv)$ is well defined for small $v$ and for small $\epsilon$ it realizes a contact isomorphism of the tangent (or cotangent) sphere of radius $\epsilon$ to the boundary of the complex tubular neighborhood of radius $\epsilon$ (cf. \cite{lB78.3}).

The corresponding FIO's can be described as follows: as above there exists a complex phase function $\phi$ on $T^*W\times T^*V^0$ such that 1) $\phi$ vanishes on the graph of $f_0$ and $d\phi=\xi.dx-\eta.dy$ there, 2) $\Im\phi\gg0$ i.e. it is positive outside of the graph and the transversal hessian is $\gg0$. $\phi$ is then  a global phase function for FIO associated to $f_0$ ($\phi$ is not unique, but obviously the set of such functions is convex, hence contractible).

The elliptic FIO's we are interested in are those that can be defined by a positive symbol (or a symbol isotopic to $1$):
$$
f\mapsto g(x)=\int e^{i\phi} a(x,\xi,y,\eta) f(y) dyd\eta d\xi  \quad\mbox{with $a>0$ on the graph}\ .
$$
The degree of such operators depends on the degree of $a$, but they all have the same index, given by the formula above.

\bk 2) The formula above extends also to vector bundle cases: if $E,E'$ are holomorphic vector bundles (or complexes of such) on $\Omega,\Omega'$, $f_0$ a contact isomorphism ($\partial\Omega\to\partial\Omega')$ as above, and $A$ a smooth (not holomorphic) isomorphism $f_{0*}E\to E'$ on the boundaries, the Toeplitz operator $a\mapsto S'(Af_{0*}a)$ is Fredholm and its index is given by the same construction as above. For this construction $f_0$ only needs to be defined where the complexes are not exact.

\mk In particular let $\Omega, \Omega'$ have singularities (isolated singularities, since we still want smooth boundaries): we can embed $\Omega, \Omega'$ in smooth strictly pseudoconvex domains $\WO,\WO'$ of the same (higher) dimension; the contact isomorphism extends at least in a small neighborhood of $\partial\Omega$ in $\partial\WO$. The coherent sheaves $\cO_\Omega,\cO_{\Omega'}$  have global locally free holomorphic resolutions on $\WO,\WO'$; near the boundary these are homotopy equivalent to a Koszul complex, hence equivalent.

The theorem above shows that the relative index is the K-theoretical character of the difference virtual bundle $d([\cO_\Omega],[\cO_{\Omega'}])$ (vanishing near the boundary). Note however that the virtual bundles $[\cO_\Omega],[\cO_{\Omega'}]$ lie in the K-theory of $\WO$ with support in $\Omega$ (resp. ..). This can be readily described in terms of cohomology classes on $\WO$ etc. with support in $\Omega$, not on $\Omega$ itself (the relation between coherent holomorphic modules and topological K-theory, or K-theory and cohomology, is not good enough when there are  singularities).

\newpage

\end{document}